\magnification=\magstep1
\input amstex
\documentstyle{amsppt}
\hoffset=.25truein
\hsize=6.0truein
\vsize=8.95truein

\topmatter
\centerline{\bf A NOTE ON SOME POSITIVITY CONDITIONS RELATED TO}
\title
zeta and L-functions
\endtitle
\keywords
Zeros of zeta and $L$-functions, Hilbert spaces of entire functions
\endkeywords
\subjclass
11M26
\endsubjclass
\author
J. B. Conrey and Xian-Jin Li
\endauthor
\address
American Institute of Mathematics, 360 Portage Avenue, Palo Alto, CA 94306
\endaddress
\thanks
Research of both authors supported by the American Institute of Mathematics.
\endthanks
\email
conrey\@aimath.org, \,\, xianjin\@math.Stanford.EDU
\endemail
\endtopmatter
\address
Current address for Xian-Jin Li: Department of Mathematics,
Brigham Young University, Provo, Utah 84602 USA
\endaddress
\document

\heading
1.   Introduction
\endheading
  The theory of Hilbert spaces of entire functions [1] was developed
by Louis de Branges in the late 1950s and early 1960s.
It is a generalization of the part of Fourier analysis involving
Fourier transforms and the Plancherel formula.  
  In [2] de Branges proposed an approach to the generalized Riemann
hypothesis, that is, the hypothesis that not only the Riemann zeta function
$\zeta(s)$ but also all the Dirichlet $L$-functions $L(s, \chi)$ with
$\chi$ primitive have their
nontrivial zeros lying on the critical line $\Re s=1/2$ (See Davenport [5]).
In [2] de Branges mentioned that his approach to the generalized Riemann
hypothesis using Hilbert spaces of entire functions is related to the
Lax-Phillips theory of scattering [6].   In Appendix 2 to Section 7, [6]
Lax and Phillips explained the difficulty of approaching the
Riemann hypothesis by using the scattering theory.
In this note, we shall indicate the difficulty of approaching
the Riemann hypothesis by using de Branges' positivity conditions [2] [3] [4].
In fact, we shall give examples showing that de Branges' positivity
conditions, which imply the generalized Riemann hypothesis, are not
satisfied by defining functions of reproducing kernel Hilbert spaces
associated with the Riemann zeta function $\zeta(s)$ and the Dirichlet
$L$-function $L(s, \chi_4)$.

\heading
2.  Reproducing kernel Hilbert spaces
\endheading

   We first outline an important part of de Branges' approach
to the Riemann hypothesis.

  Let $E(z)$ be an entire function satisfying $|E(\bar z)|<|E(z)|$ for
$z$ in the upper half-plane.  A Hilbert space of entire functions
$\Cal H(E)$ is the set of all entire functions $F(z)$ such that
$F(z)/E(z)$ is square integrable on the real axis and such that
$$|F(z)|^2\leqslant \|F\|^2_{\Cal H(E)} K(z, z) \tag 2.1$$
for all complex $z$, where the inner product of the space is
given by
$$\langle F(z), G(z)\rangle_{\Cal H(E)}=\int_{-\infty}^\infty
{F(x)\bar G(x)\over |E(x)|^2}dx$$
for all elements $F, G\in\Cal H(E)$ and where
$$K(w, z)={E(z)\bar E(w)-\bar E(\bar z)E(\bar w)\over
2\pi i(\bar w-z)}$$
is the reproducing kernel function of the space $\Cal H(E)$,
that is, the identity
$$F(w)=\langle F(z), K(w, z)\rangle_{\Cal H(E)} \tag 2.2$$
holds for every complex $w$ and for every element $F\in\Cal H(E)$.
The identity (2.2) is obtained by using Cauchy's integration
formula in the upper half-plane (cf. [1]), and the condition
(2.1) is made so that Cauchy's formula applies to all functions
in the space $\Cal H(E)$.

   The following two theorems are essentially due to de Branges (cf. [2] [3]).

 \proclaim{Theorem 1}   Let $E(z)$ be an entire function having no real zeros
such that $|E(\bar z)|<|E(z)|$ for $\Im z>0$, such that $\bar E(\bar z)=\epsilon E(z-i)$
for a constant $\epsilon$ of absolute value one, and such that
$|E(x+iy)|$ is a strictly increasing function of $y>0$ for each fixed real $x$.
If  $\Re\langle F(z), F(z+i)\rangle_{\Cal H(E)}\geqslant 0$ for every
element $F(z)\in\Cal H(E)$ with $F(z+i)\in\Cal H(E)$, then the zeros
of $E(z)$ lie on the line $\Im z=-1/2$, and
$\Re\{\bar{E^\prime}(w) E(w+i)/2\pi i\}\geqslant 0$
when $w$ is a zero of $E(z)$. \endproclaim

\demo{Proof}  Let $w$ be a zero of $E(z)$.  Since
$\bar E(\bar z)=\epsilon E(z-i)$ with $|\epsilon|=1$,
we have
$$\bar E(w+i)K(w, z+i)=-\bar E(w-i)K(w+i, z)\tag 2.3$$
for all complex $z$.  Since $E(z)$ has no real zeros
and since $|E(\bar z)|<|E(z)|$ for $z$ in the upper half-plane,
$E(w+i)$ and $E(w-i)$ are nonzero.  It follows that $K(w, z)$ is a
nonzero element of $\Cal H(E)$ such that $K(w, z+i)$ belongs to the space.
Assume that $F(z)$ is an element in $\Cal H(E)$ such that $F(z+i)$
belongs to the space.  Then, by (2.2) and (2.3), we have
$$\aligned &\langle F(z+i), K(w, z)\rangle_{\Cal H(E)}
+\langle F(z), K(w, z+i)\rangle_{\Cal H(E)}\\
&= {E(w+i)-E(w-i)\over E(w+i)}F(w+i). \endaligned\tag 2.4$$
Define a new scalar product $\langle \cdot, \cdot\rangle$ by
$$\langle F(z), G(z)\rangle=\langle F(z+i), G(z)\rangle_{\Cal H(E)}
+\langle F(z), G(z+i)\rangle_{\Cal H(E)}$$
for all $F, G\in\Cal H(E)$ such that $F(z+i), G(z+i)\in\Cal H(E)$.
Since, by assumption, $\Re\langle F(z+i), F(z)\rangle_{\Cal H(E)}\geqslant 0$
for every element $F(z)\in\Cal H(E)$ such that $F(z+i)\in\Cal H(E)$, we have
$\langle F(z), F(z)\rangle\geqslant 0$.   Then, by (2.4) and the Schwarz
inequality, we have
$$\aligned &\left|{E(w+i)-E(w-i)\over E(w+i)}F(w+i)\right|^2
=|\langle F(z), K(w, z)\rangle|^2\\
&\leqslant \langle F(z), F(z)\rangle \langle K(w,z), K(w, z)\rangle\\
&=4Re\, K(w, w+i)Re\,\langle F(z+i), F(z)\rangle_{\Cal H(E)}. \endaligned \tag 2.5$$

   If $K(w, w+i)\neq 0$, then we must have $\Im w=-1/2$ because, otherwise,
we have $K(w, w+i)=0$ by the functional identity $\bar E(\bar z)=\epsilon E(z-i)$.

   Next, we assume that $K(w, w+i)=0$.  Then $F(z)=K(w, z)$ is an element of
$\Cal H(E)$ such that $F(w+i)=0$ and $F(z+i)\in\Cal H(E)$.
If $F(z)$ is a nonzero element of $\Cal H(E)$
having zero at a point $z_0$, it is easy to see by definition that $F(z)/(z-z_0)$
belongs to $\Cal H(E)$.  Since $F(z)$ is an entire function, by using
Taylor's expansion of $F(z)$ at the point $z_0$, we see that $F(z)/(z-z_0)^n$
does not vanish at $z_0$ for some positive integer $n$.  By the repeated process
of dividing out the factor $z-z_0$ from $F(z)$, we see that $F(z)/(z-z_0)^n$
belongs to $\Cal H(E)$.   If $F(z+i)\in\Cal H(E)$, we also see that
$F(z+i)/(z+i-z_0)^n\in\Cal H(E)$.
 Therefore, there exists a nonzero element $F(z)\in\Cal H(E)$
such that $F(w+i)\neq 0$ and $F(z+i)\in\Cal H(E)$.
Hence, if $K(w, w+i)=0$, then, by (2.5), we have
$$E(w+i)-E(w-i)=0.$$
Since $\bar E(\bar w)=\epsilon E(w-i)$, we have $|E(w+i)|=|E(\bar w)|$.
Note that $\Re(w+i)=\Re(\bar w)$.  Since $|E(x+iy)$ is a strictly
increasing function of $y$ on $(0, \infty)$, we must have
$\Im(w+i)=\Im(\bar w)$, and hence $w+i=\bar w$.
Therefore, we have $\Im w= -{1\over 2}$.

 We have $K(w, w+i)=\bar{E^\prime}(w)E(w+i)/2\pi i$.
Since $F(z)=K(w, z)$ is an element of $\Cal H(E)$ such that
$F(z+i)\in\Cal H(E)$, we have
$$\Re K(w, w+i)=\Re\langle F(z+i), F(z)\rangle_{\Cal H(E)}\geqslant 0,$$
that is, $\Re\{\bar{E^\prime}(w) E(w+i)/2\pi i\}\geqslant 0$
when $w$ is a zero of $E(z)$.

 This completes the proof of the theorem.\enddemo

   Let $W(z)$ be a function analytic and having no zeros in
the upper half-plane.  Then a Hilbert space of analytic functions $\Cal F(W)$
is the set of all analytic functions $F(z)$ in the upper half-plane, such that
$F(z)/W(z)$ can be written as a quotient of bounded analytic functions in the
upper half-plane, has square integrable boundary values on the real axis, and
satisfies the inequality
$$\log|F(x+iy)/W(x+iy)|\leqslant \frac{y}\pi\int^{+\infty}_{-\infty}
\frac{\log|F(t)/W(t)|dt}{(t-x)^2+y^2}$$
for $y>0$.  The inner product of $\Cal F(W)$ is given by
$$\langle F(z), G(z)\rangle_{\Cal F(W)}=\int_{-\infty}^\infty
{F(x)\bar G(x)\over |W(x)|^2}dx$$
for all $F, G\in\Cal F(W)$.  The reproducing kernel function of $\Cal F(W)$ is
given by the expression
$$K(w, z)=\frac{W(z)\bar W(w)}{2\pi i(\bar w-z)},$$
that is, for every complex $w$ in the upper half-plane, we have
$$F(w)=\langle F(z), K(w, z)\rangle_{\Cal F(W)} \tag 2.6$$
for every element $F\in\Cal F(W)$.
The identity (2.6) is obtained by using Cauchy's integration
formula in the upper half-plane (cf. [1]).

\proclaim{Theorem 2}  Let $W(z)$ be a function analytic and having no
zeros in the upper half-plane.   Let $T$ be a linear transformation of
$\Cal F(W)$ into itself which takes $K(w, z)$ into $K(w+i, z)$
for all complex $w$ with $\Im w>0$.  Assume that
$$\Re\langle F(z), TF(z)\rangle_{\Cal F(W)}\geqslant 0$$
for all $F\in\Cal F(W)$.  Then $W(z)$ has an analytic extension to
the half-plane $\Im z>-1/2$, and $W(z)/W(z+i)$ has a nonnegative
real part in this half-plane.
\endproclaim

\demo{Proof}   Let $w_1,\dots,w_r$ be points in the upper half-plane, and
let $c_1,\dots,c_r$ be complex numbers.   If
$F(z)=\sum_{\alpha=1}^rc_\alpha K(w_\alpha, z)$, then
$TF(z)=\sum_{\beta=1}^rc_\beta K(w_\beta +i, z)$.
By assumption, we have
 $$\sum_{\alpha, \beta=1}^r c_\alpha\bar c_\beta
[K(w_\alpha, w_\beta +i)+K(w_\alpha+i, w_\beta)]
=2\Re\langle F(z), TF(z)\rangle_{\Cal F(W)}\geqslant 0, $$
that is, the expression $K(w+i, z)+K(w, z+i)$ is positive-definite
for $w, z$ in the upper half-plane.  This implies that
$\Re \{W(z)/W(z+i)\}\geqslant 0$ for $z$ in the upper half-plane.
Let
$$B(z)={W(z)-W(z+i)\over W(z)+W(z+i)}.$$
Then $B(z)$ is analytic and bounded by one in the upper half-plane.
The positive-definiteness of $K(w+i, z)+K(w, z+i)$
implies the positive-definiteness of the expression
$${1-B(z)\bar B(w)\over 2\pi i(\bar w-z-i)}$$
for $w, z$ in the upper half-plane.

 Let $\Cal H$ be the Hilbert space of analytic functions in the half-plane
$\Im z>-1/2$, which has the expression $L(w, z)=1/2\pi i(\bar w-z-i)$ as its
 reproducing kernel function.  The norm of an
element $F$ in the space $\Cal H$ is given by
$$\langle F(z), G(z)\rangle_{\Cal H}=\int_{-\infty}^\infty
F(x-i/2)\bar G(x-i/2)dx.$$

  Let $P$ be a transformation of $\Cal H$ into itself, which takes
$L(w, z)$ into $\bar B(w)L(w, z)$.
 The positive-definiteness of the expression
$${1-B(z)\bar B(w)\over 2\pi i(\bar w-z-i)}$$
implies $\langle PF(z), PF(z)\rangle_{\Cal H}\leqslant \langle F(z),
F(z)\rangle_{\Cal H}$ for all elements $F\in\Cal H$ which are linear
combination of functions $L(w, z)$ with $\Im w>0$.  Since $L(w, z)$
is the reproducing kernel function of $\Cal H$, if $F\in\Cal H$ is
orthogonal to all elements $L(w, z)$ with $\Im w>0$, then
$$F(w)=\langle F(z), L(w, z)\rangle_{\Cal H}=0$$
for $\Im w>0$.  Since $F(z)$ is analytic for $\Im z>-1/2$,
we must have $F\equiv 0$.  Therefore, the set of elements $L(w, z)$
with $\Im w>0$ is dense in $\Cal H$.  It follows that
$\langle PF(z), PF(z)\rangle_{\Cal H}\leqslant \langle F(z),
F(z)\rangle_{\Cal H}$ for all elements $F\in\Cal H$.  Thus,  $P$ is
a bounded linear transformation of the Hilbert space $\Cal H$
into itself, and therefore, the adjoint $P^*$ of $P$ exists.

    Let $\alpha$ be a complex number with $\Im\alpha>-1/2$,
and let $F(z)=L(\alpha, z)$.  Then $F\in \Cal H$, and hence
$P^*F(z)\in\Cal H$.  It follows that
$$2\pi i(\bar\alpha-w-i)\langle P^*F(z), L(w, z)\rangle_{\Cal H}
=\int_{-\infty}^\infty P^*F(x-i/2){\bar\alpha-w-i\over
x-{i\over 2}-w}dx$$
is an analytic function of $w$ for $\Im w>-1/2$.  Since
$$2\pi i(\bar\alpha-w-i)\langle P^*F(z), L(w, z)\rangle_{\Cal H}
=2\pi i(\bar\alpha-w-i)\langle F(z), PL(w, z)\rangle_{\Cal H}=B(w)$$
for $w$ in the upper half-plane, $B(z)$ has an analytic extension to the
half-plane $\Im z>-1/2$.

  If $F\in\Cal H$ and $\Im w>0$, we have
$$B(w)F(w)=\langle P^*F(z), L(w, z)\rangle_{\Cal H}. \tag 2.7$$
Since both sides of (2.7) are analytic functions of $w$ for
$\Im w>-1/2$, the identity (2.7) remains true for all complex $w$
with $\Im w>-1/2$ by analytic continuation.  Since
$\langle PF(z), PF(z)\rangle_{\Cal H}\leqslant \langle F(z),
F(z)\rangle_{\Cal H}$ for all $F\in\Cal H$, we have
$$\aligned |B(w)F(w)|^2&=|\langle F(z), PL(w, z)\rangle_{\Cal H}|^2\\
&\leqslant \langle F(z), F(z)\rangle_{\Cal H}\langle PL(w, z),
PL(w, z)\rangle_{\Cal H}\\
&\leqslant \langle F(z), F(z)\rangle_{\Cal H}\langle L(w, z),
L(w, z)\rangle_{\Cal H}\\
&=\langle F(z), F(z)\rangle_{\Cal H}L(w, w).\endaligned \tag 2.8$$
for $\Im w>-1/2$ and for all $F\in\Cal H$.  In particular, if
$F(z)=L(w,z)$, then $\langle F(z), F(z)\rangle_{\Cal H}=F(w)$, and hence
(2.8) becomes $|B(w)F(w)|^2\leqslant |F(w)|^2$ for
$\Im w>-1/2$, that is, $|B(w)|\leqslant 1$ for $\Im w>-1/2$.

  Therefore, we have proved that $B(z)$ is analytic and
bounded by one for $\Im z>-1/2$.  It follows that $W(z)/W(z+i)$ is analytic
and has nonnegative real part in the half-plane $\Im z>-1/2$.

This completes the proof of the theorem.
\enddemo

\heading
3.  Hilbert spaces associated with $\zeta(s)$ and $L(s, \chi_4)$
\endheading

\subheading{3.1.  The Riemann zeta function}   The Riemann zeta function $\zeta(s)$
is given by
$$\zeta(s)=\sum_{n=1}^\infty {1\over n^s}$$
for $\Re s>1$.    Let $\xi(s)=s(s-1)\pi^{-s/2}\Gamma({s\over 2})\zeta(s)$.
Then $\xi(s)$ is an entire function, and satisfies the functional identity
$\xi(s)=\xi(1-s)$.   It is well-known (See Davenport [5]) that we have
the infinite product formula
$$\xi(s)=\prod \left(1-{s\over\rho}\right)$$
where the product is taken over all nontrivial zeros $\rho$ of $\zeta(s)$
with $\rho$ and $1-\rho$ being paired together for the convergence of
the product.

    Let $E(z)=\xi(1-iz)$.  Then the Riemann hypothesis is that the zeros
of $E(z)$ lie on the line $\Im z=-1/2$, and the functional identity
$\xi(s)=\xi(1-s)$ can be written as $\bar E(\bar z)=E(z-i)$.  If
$\rho$ is a nontrivial zero of $\zeta(s)$, then $0<\Re\rho<1$, a
result proved by Hadamard and de la Vall\'ee Poussin independently in 1896.   Since
$$|E(z)|^2=\prod \left|1-{iz\over\rho}\right|^2
=\prod {(\Re\rho+y)^2+(\Im\rho-x)^2 \over |\rho|^2}$$
for $z=x+iy$, we see that $|E(x-iy)|<|E(x+iy)|$ for $y>0$, and that
$|E(x+iy)|$ is a strictly increasing function of $y$ on $(0, \infty)$ for
each fixed real $x$.

    In view of Theorem 1, it is natural to ask whether the Hilbert
space of entire functions $\Cal H(E)$ satisfies the condition that
$$\Re \langle F(z), F(z+i)\rangle_{\Cal H(E)}\geqslant 0 \tag 3.1$$
for every element $F(z)$ of $\Cal H(E)$ such that $F(z+i)\in\Cal H(E)$,
because the nontrivial zeros of the Riemann zeta function $\zeta(s)$
would then lie on the critical line $\Re s=1/2$ under this condition.
In the following, we give an example showing that condition (3.1) is
unfortunately not true.

   Let $\rho=1/2+i111.0295355431696745\cdots$ be the $34$th zero of the Riemann
zeta function in the upper half-plane.  By using {\it MATHEMATICA},
we compute that
$$-\Re\{\xi^\prime(\rho)\xi(1+\rho)\}=-5.389100507182945\cdots\times 10^{-69}
<0. \tag 3.2$$
Write $\rho=1-iw$.  Then $E(w)=0$, and $\bar {E^\prime}(w)E(w+i)/i
=-\xi^\prime(\rho)\xi(1+\rho)$.  Thus, (3.2) becomes
$$\Re\{\bar {E^\prime}(w)E(w+i)/2\pi i\}<0.$$
Therefore, by Theorem 1, we see that the Hilbert space of entire
functions $\Cal H(E)$ with $E(z)=\xi(1-iz)$ does not satisfy
the condition (3.1).

   Next, let $W(z)=1/\xi(1-iz)$.   Then $W(z)$ is analytic in the upper
half-plane, and is continuous
and having no zeros in the closed upper half-plane.  If the Hilbert
space of analytic functions $\Cal F(W)$ satisfy the condition that
$$\Re\langle F(z), TF(z)\rangle_{\Cal F(W)}\geqslant 0 \tag 3.3$$
for all $F\in\Cal F(W)$, where $T$ is the linear transformation of
$\Cal F(W)$ into itself which takes $K(w, z)=W(z)\bar W(w)/2\pi i(\bar w-z)$
into $K(w+i, z)$ for all complex $w$ with $\Im w>0$, then the function
$W(z)$ would have analytic extension to the half-plane $\Im z>-1/2$
by Theorem 2, that is, the Riemann zeta function $\zeta(s)$ has no
zeros for $\Re s>1/2$.   In the following, we give an example showing that
the space $\Cal F(W)$ does not satisfy the condition (3.3).

  By using {\it MATHEMATICA}, we compute that
$$\Re\{\xi(1+i282)/\xi(2+i282)\}=-0.000131957<0. \tag 3.4$$
Let $w=-282$.  Then $\Im w=0>-1/2$.  Since
 $W(w)=1/\xi(1+i282)$ and $W(w+i)=1/\xi(2+i282)$, by (3.4) we have
$$\Re\{W(w)/W(w+i)\}<0.$$
Therefore, by Theorem 2, we see that the space
$\Cal F(W)$ with $W(z)=1/\xi(1-iz)$
 does not satisfy the condition (3.3).

  \subheading{3.2.  The Dirichlet $L$-function $L(s, \chi_4)$}
Note that $\chi_4$ is the real primitive Dirichlet character
(mod $4$), which is given by
$$\chi_4(n)=\cases (-1)^{n-1\over 2}, &\text{if $n$ is odd;}\\
0, &\text{if $n$ is even}.\endcases$$
The Dirichlet $L$-function $L(s, \chi_4)$ is given by
$$L(s, \chi_4)=\sum_{n=1}^\infty {\chi_4(n)\over n^s}$$
for $\Re s>0$.  Let $\xi(s, \chi_4)=(4/\pi)^{s\over 2}\Gamma({1+s\over 2})
L(s, \chi_4)$.  Then $\xi(s, \chi_4)$ is an entire function, and satisfies
the functional identity
$$\xi(1-s, \chi_4)=\epsilon(\chi_4)\xi(s, \chi_4) \tag 3.5$$
where $\epsilon(\chi_4)$ is a constant of absolute value one
(See Davenport [5]).   Since $\chi_4$ is a real character, by the
argument of section 12 of Davenport [5], we have the infinite product
formula
$$\xi(s, \chi_4)={\sqrt\pi\over 2}\prod \left(1-{s\over\rho}\right) \tag 3.6$$
where  the product is taken over all nontrivial
zeros of $L(s, \chi_4)$ with $\rho$ and $1-\rho$ being put together.

   Let $E_{\chi_4}(z)=\xi(1-iz, \chi_4)$.  Then the functional identity
(3.5) can be written as
$$\bar E_{\chi_4}(\bar z)=\bar\epsilon(\chi_4)E_{\chi_4}(z-i).$$
  By the product formula (3.6), we find that
$|E_{\chi_4}(x+iy)|$ is a strictly increasing function of $y>0$
for each fixed real $x$.  Since the nontrivial zeros of $L(s, \chi_4)$
lie in the strip $0<\Re s<1$ (See \S14, Davenport [5]), we
have $|E_{\chi_4}(\bar z)|<|E_{\chi_4}(z)|$ for $\Im z>0$.

    In view of Theorem 1, it is natural to ask whether the Hilbert
space of entire functions $\Cal H(E_{\chi_4})$ satisfies the condition that
$$\Re \langle F(z), F(z+i)\rangle_{\Cal H(E_{\chi_4})}\geqslant 0 \tag 3.7$$
for every element $F(z)$ of $\Cal H(E_{\chi_4})$ such that
$F(z+i)\in\Cal H(E_{\chi_4})$, because the nontrivial zeros of
the Dirichlet $L$-function $L(s, \chi_4)$
would then lie on the critical line $\Re s=1/2$ under this condition.
In the following, we give an example showing that condition (3.7) is
unfortunately not true.

   Let $\rho=1/2+i67.6369208635460683980549\cdots$ be a zero of $L(s, \chi_4)$.
 By using {\it MATHEMATICA}, we compute that
$$-\Re\{\xi^\prime(\rho, \chi_4)\xi(1+\rho, \chi_4)\}
=-2.310349004993483456\cdots\times 10^{-45}<0.$$
Write $\rho=1-iw$.  Then $E_{\chi_4}(w)=0$, and
$\bar {E_{\chi_4}^\prime}(w)E_{\chi_4}(w+i)/i
=-\xi^\prime(\rho, \chi_4)\xi(1+\rho, \chi_4)$.  Thus,
the above inequality becomes
$$\Re\{\bar {E_{\chi_4}^\prime}(w)E_{\chi_4}(w+i)/2\pi i\}<0.$$
Therefore, by Theorem 1, we see that the Hilbert space of entire
functions $\Cal H(E_{\chi_4})$ with $E_{\chi_4}(z)=\xi(1-iz, \chi_4)$
does not satisfy the condition (3.7).

  Next, let $W_{\chi_4}(z)=1/\xi(1-iz, \chi_4)$.   Then $W_{\chi_4}(z)$
is analytic in the upper half-plane, and is continuous
and having no zeros in the closed upper half-plane.  If the space
 $\Cal F(W_{\chi_4})$ satisfy the condition that
$$\Re\langle F(z), TF(z)\rangle_{\Cal F(W_{\chi_4})}\geqslant 0 \tag 3.8$$
for all $F\in\Cal F(W_{\chi_4})$, where $T$ is the linear transformation of
$\Cal F(W_{\chi_4})$ into itself which takes $K(w, z)=W_{\chi_4}(z)
\bar W_{\chi_4}(w)/2\pi i(\bar w-z)$
into $K(w+i, z)$ for all complex $w$ with $\Im w>0$, then the function
$W_{\chi_4}(z)$ would have analytic extension to the half-plane $\Im z>-1/2$
by Theorem 2, that is, the Dirichlet $L$-function $L(s, \chi_4)$ has no
zeros for $\Re s>1/2$.   In the following, we give an example showing that
the space $\Cal F(W_{\chi_4})$ does not satisfy the condition (3.8).

   By using {\it MATHEMATICA}, we compute that
$$\Re\{\xi(1+i8714.2, \chi_4)/\xi(2+i8714.2, \chi_4)\}=-0.000422340607<0. $$
Let $w=-8714.2$.  Then $\Im w=0>-1/2$.  Since
 $W_{\chi_4}(w)=1/\xi(1+i8714.2, \chi_4)$ and
$W_{\chi_4}(w+i)=1/\xi(2+i8714.2, \chi_4)$, we have
$$\Re\{W_{\chi_4}(w)/W_{\chi_4}(w+i)\}<0. $$
Therefore, by Theorem 2, we see that the space
$\Cal F(W_{\chi_4})$ with $W_{\chi_4}(z)=1/\xi(1-iz, \chi_4)$
 does not satisfy the condition (3.8).

\heading
4.  Conclusion
\endheading

   We have seen in section 3 the difficulty of approaching the
generalized Riemann hypothesis by using de Branges type
positivity conditions for reproducing kernel Hilbert spaces.
It is possible that these positivity conditions are too strong for
Hilbert spaces of entire functions associated with the Riemann zeta
function and the Dirichlet $L$-functions.

\vskip0.25truein
  \subheading{Remark}    After he looked at the manuscript
of this paper, {\it Peter Sarnak} gave a proof for the statement 
that the space $\Cal F(W)$ does not satisfy the condition (3.3) 
where $W(z)=1/\xi(1-iz)$, and his argument involves no
numberical calculations.   For the convenience of
readers, we sketch his proof here.
 Let $F(s)=\xi(s)/\xi(s+1)$.  Assume that all logarithmic 
functions are defined by continuation from the point
$s=2$ at which their arguments are set to be zero.
Then we have
$$\Im \{\log F(s)\}=\Im\{\log \zeta(s)\}+O(1)$$
for $\Re s>1/2$.   Since the set of values of
$\log \zeta(s)$, $1/2 <\Re s < 2$, is dense in the complex plane
(See Chapter XI of Titchmarsh [7]),  a complex number
 $s_0$ exists with $\Re s_0>1/2$ such that 
$\pi/2<\Im\{\log F(s_0)\}<\pi$.  Let $z_0=i(s_0-1)$.  
Then $\Im z_0>-1/2$ and $\Re\{W(z_0)/W(z_0+i)\}< 0$, 
and hence by Theorem 2 the space $\Cal F(W)$ 
 does not satisfy the condition (3.3).
Let $r$ be any positive integer.  For any Dirichlet character 
$\chi$ modulo $r$, let
$$\xi(s, \chi)=(\pi/r)^{-{s+a\over 2}}
\Gamma({s+a\over 2})L(s, \chi)$$
where $L(s, \chi)$ is the Dirichlet $L$-function
and where $a=0$ if $\chi(-1)=1$ and $a=1$ if $\chi(-1)=-1$.
  By using a similar argument,
{\it Peter Sarnak} also proved that the space $\Cal F(W_\chi)$ 
does not satisfy the condition (3.8) 
where $W_\chi(z)=1/\xi(1-iz, \chi)$.

\vskip0.15truein
\heading
 Appendix
\endheading

  Here is a list of {\it MATHEMATICA} instructions, which can be used to
verify our calculations.

 {\tt
xi[s\_]:=s*(s-1)*Pi\^\,(-s/2)*Gamma[s/2]*Zeta[s]

f[ro\_]:=Re[-xi'[ro]*xi[ro+1]]

g[t\_]:=Re[xi[1+I*t]/xi[2+I*t]]

ro1=1/2+I*111.029535543169674524656

Zeta[ro1]

f[ro1]

g[282.]

Plot[g[t],\{t,281.95,282.15\}]

xi4[s\_]:=(4*Pi)\^\,(-s/2)*Gamma[(s+1)/2]*(Zeta[s,1/4]-Zeta[s,3/4])

f4[ro\_]:=Re[-xi4'[ro]*xi4[ro+1]]

g4[t\_]:=Re[xi4[1+I*t]/xi4[2+I*t]]

ro2=1/2+I*67.6369208635460683980549

Zeta[ro2,1/4]-Zeta[ro2,3/4]

f4[ro2]

g4[8714.2]

Plot[g4[t],\{t,8714.1,8714.4\}]
}
\enddocument